\documentclass[12pt]{article}

\usepackage{pict2e}
\usepackage{amsmath}
\usepackage{amssymb}
\author{Abdallah Assi\thanks{Universit\'e d'Angers, D\'epartement de
Math\'ematiques, 2 bd Lavoisier, 49045 Angers Cedex 01, France, e-mail:assi@univ-angers.fr}}
\title{On subsemigroups of ${\mathbb{N}}^e$}
\date{\mbox{}}

\newtheorem{teorema}{Theorem}[section]

\newtheorem{proposicion}[teorema]{Proposition}

\newtheorem{definicion}[teorema]{Definition}

\newtheorem{corolario}[teorema]{Corollary}

\newtheorem{nota}[teorema]{Remark}
\newtheorem{exemple}[teorema]{Example}

\newenvironment{demostracion}[1]{\paragraph{\sl Proof#1}}{}

\begin{document}
\maketitle

\begin{center}
{\bf Abstract}
\end{center}
\noindent{\small Let $\underline{v}=(v_1,\cdots,v_{e+s})$ be a set of vectors of ${\mathbb{N}}^e$, and assume that
 $v_{e+k}$ is not in the group generated by $v_1,\ldots,v_{e+k-1}$ for all $k=1,\cdots,s$. The aim of this paper is
 to  give a formula for the Frobenius number and the
conductor of the subsemigroup generated par $\underline{v}$ in $\mathbb{N}^e$.}

\section{Introduction and Basic Notations}

\noindent Let $\underline{v}=(v_1,\ldots,v_e,v_{e+1},\ldots,v_{e+s})$ be a set of
nonzero elements of $\mathbb{N}^e$ and let

$$
\Gamma(\underline{v})=\lbrace
\sum_{k=1}^{e+s}a_iv_i|a_i\in\mathbb{N}\rbrace
$$

\noindent be the subsemigroup of
$\mathbb{N}^e$ generated by $\underline{v}$. Let $G(\underline{v})=\lbrace
\sum_{k=1}^{e+s}a_iv_i|a_i\in\mathbb{Z}\rbrace$ be the subgroup of
${\mathbb Z}^e$ generated
by $\underline{v}$ and let cone$(v_1,\ldots,v_e)$ be the convex cone
generated by $v_1,\ldots,v_e$,

$$
{\rm cone}(v_1,\cdots,v_e)=\lbrace
\sum_{k=1}^{e}a_iv_i|a_i\in\mathbb{R}_+\rbrace.
$$

\noindent Assume that the
dimension of cone$(v_1,\ldots,v_e)$ is $e$ -i.e. $(v_1,\ldots,v_e)$
generates $\mathbb{R}^e$- and that $v_{e+1},\ldots,v_{e+s}\in {\rm cone}(v_1,\ldots,v_e)$. The
paper deals with the following question:

\vskip0.1cm

\noindent What is the ``smallest'' element $w\in {\rm cone}(v_1,\ldots,v_e)$ such that for
all $v\in w+{\rm cone}(v_1,\ldots,v_e)$, if $v\in G(\underline{v})$,
then $v\in \Gamma(\underline{v})$?

\vskip0.1cm

\noindent Let $D_1$ be the determinant of the matrix
$[v_1^T,\ldots,v_e^T]$ -where $T$ denotes the transpose of a matrix-,
and for all $k=2,\ldots,s+1$, let $D_k$ be the gcd of the $(e,e)$ minors
of the matrix $[v_1^T,\ldots,v_e^T,v_{e+1}^T,\ldots,v_{e+k-1}^T]$. Set $e_k=\dfrac{D_k}{D_{k+1}}$ for all $k=1,\ldots,s$. We shall assume that the two following conditions are satisfied:

(*) $D_1 > D_2 >\cdots >D_{s+1}$ (in particular for
all $k=2,\cdots,s+1, v_{e+k-1}$ is not in the group generated by
$(v_1,\cdots,v_e,v_{e+1},\cdots,v_{e+k-2})$).

(**)  $e_kv_{e+k}\in \Gamma(v_1,\ldots,v_e,v_{e+1},\ldots,v_{e+k-1})$ for all $k=1,\ldots,s$.

\noindent Our main result is
the following:

\vskip0.1cm

\noindent{\bf Theorem 1.} (See Fig. 1) Let the notations be as above, and let $C_e$
be the unique cell of dimension $e$ of cone$(v_1,\ldots,v_e)$ (more
precisely $C_e$ is the interior of cone$(v_1,\ldots,v_e)$). If

$$
g(\underline{v})=\sum_{k=1}^s (e_k-1)v_{e+k}-\sum_{i=1}^e v_i
$$

\noindent then the following hold:

\vskip0.1cm

i) $g(\underline{v})\notin \Gamma(\underline{v})$.

ii) For
all $v\in g+(C_e-\lbrace (0,\cdots,0)\rbrace)$, if $v\in G(\underline{v})$,
then $v\in \Gamma(\underline{v})$.

\noindent We call $g(\underline{v})$ the Frobenius vector of $\Gamma(\underline{v})$.

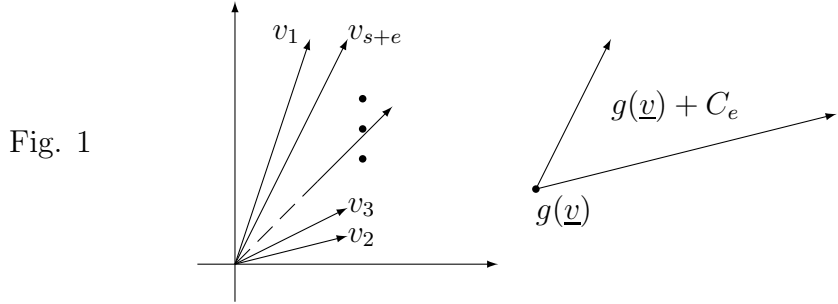
 \unitlength=1cm
\begin{picture}(6,4)

\put(6,-0.5){\vector(0,1){4}} 
\put(6.5,3){\shortstack{$v_1$}}
\put(7.5,0.3){\shortstack{$v_2$}}
\put(7.5,0.7){\shortstack{$v_3$}}
\put(7.7,1.4){\circle*{.1}}
\put(7.7,1.8){\circle*{.1}}
\put(7.7,2.2){\circle*{.1}}
\put(7.5,3){\shortstack{$v_{s+e}$}}
\put(6,0){\line(1,1){0.2}}
\put(6.3,0.3){\line(1,1){0.2}}
\put(6.6,0.6){\line(1,1){0.2}}
\put(6.9,0.9){\line(1,1){0.2}}
\put(7.1,1.1){\vector(1,1){1}}
\put(6,0){\vector(1,2){1.5}}  
\put(6,0){\vector(2,1){1.5}}  
\put(6,0){\vector(1,3){1}}
\put(6,0){\vector(4,1){1.5}}
\put(10,1){\vector(4,1){4}}
\put(10,1){\vector(1,2){1}}
\put(10,0.6){\shortstack{$g(\underline{v})$}}
\put(11,2){\shortstack{$g(\underline{v})+C_e$}}
\put(10,1){\circle*{.1}}
\put(3,1.5){\shortstack{Fig. 1}}
\put(5.5,0){\vector(1,0){4}}
\end{picture}

\bigskip

\noindent When $e=s=1$ and $v_1,v_2$ are relatively prime elements of ${\mathbb N}$, Sylvester proved in [9] that the Frobenius number of $\Gamma(v_1,v_2)$
is $(v_1-1).(v_2-1)-1$ (note that in this case, $e_1=v_1$). In [6], M.J. Knight generalized the formula for Sylvester
to a system of elements $(v_1,\cdots,v_e,v_{e+1})\in{\mathbb  N}^e$, assuming that $v_{e+1}\in{\rm cone}(v_1,\cdots,v_e)$, and that
$v_1,\cdots,v_e,v_{e+1}$ generate ${\mathbb Z}^e$. Hence, Theorem 1. can be viewed as  a generalisation of Knight's result.

\medskip

\noindent Let $e=1$ and assume that $v_1 <\cdots < v_{s+1}$. Assume, without loos of generality, that $D_{s+1}=1$, i.e. $v_1,\cdots,v_{s+1}$ are relatively prime.
The above theorem says that for all $v\geq g+1=\sum_{k=1}^s (e_k-1)v_{1+k}-v_1+1, v\in \Gamma(\underline{v})$. The
positive integer $g+1$ is called the conductor of
$\Gamma(\underline{v})$ in $\mathbb{N}$. In fact, the ideal $(t^c)$ is
the conductor ideal of the algebra ${\bf K}[t^{v_1},\cdots,t^{v_{s+1}}]$
over a field ${\bf K}$
into its integral closure ${\bf K}[t]$.

\medskip

\vskip0.1cm

\section{Lattices in $\mathbb{Z}^e$}

\noindent Let the notations be as in Section 1. The group $G_k=G(v_1,\cdots,v_e,v_{e+1},\cdots,v_{e+k}))$ being
a subgroup of the free group $\mathbb{Z}^e$ for all $k=0,\cdots,s$, it follows that $G_k$ is
a free group of rank $\leq e$, and the hypothesis on $v_1,\ldots,v_e$ implies that
the rank of $G_k$ is $e$. Let
$w_1,\ldots,w_e$ be a basis of $G_k$, in particular $D_{k+1}$ is the
determinant of the $(e,e)$ matrix $[w_1^T,\ldots,w_e^T]$. Furthermore we have the
following:

\begin{proposicion}{\rm  Let $v$ be a nonzero element of $\mathbb{Z}^e$ and
  denote by $D$ the gcd of the $(e,e)$ minors of the matrix
$[v_1^T,\ldots,v_{e+k}^T,v^T]$. Then $D$ is also the gcd of the $(e,e)$
minors of the matrix  $[w_1^T,\ldots,w_e^T,v^T]$. We also have the following:

\vskip0.2cm

 i) $D$ divides $D_{k+1}$ and $v\in G_k$ if and only if $D_{k+1}=D$.
\vskip0.2cm
 ii) $\displaystyle{{D_{k+1}\over D}}.v\in G_k$ and if $D_{k+1} > D$ then
for all $1\leq i<\displaystyle{{D_{k+1}\over D}}, i.v\notin G_k$.}
\end{proposicion}

\begin{demostracion}{.} i) For all $i=1,\ldots,e$, let $d_i$ be the determinant of the matrix
  $[w_1^T,\ldots,w_{i-1}^T,v^T,w_{i+1}^T,\ldots,w_e^T]$ and note that $D$ divides $d_i$. If $D=D_{k+1}$ then $D_{k+1}$
divides $d_i$ for all $1\leq i\leq e$. In particular the system
  $\lambda_1w_1+\ldots+\lambda_ew_e=v$ has the unique solution
  $\lambda_i=\displaystyle{{d_i\over D_{k+1}}}\in {\mathbb{Z}}$. Conversely, if $v\in G_k$, then
  there exist unique integers $\lambda_1,\ldots,\lambda_e$ such that
  $v=\lambda_1w_1+\ldots+\lambda_ew_e$, but
  $(\lambda_1,\ldots,\lambda_e)$ is the unique solution to the $e\times e$ system
  $a_1w_1+\ldots+a_ew_e=v$, in particular $\lambda_i=\displaystyle{{d_i\over D_{k+1}}}$, and $D_{k+1}$ divides $d_i$
  for all $i=1,\ldots,e$. Since $D=$ gcd$(d_1,\cdots,d_e,D_{k+1})$, then $D=D_{k+1}$.

ii) Let the notations be as in i) and $1\leq i
<\displaystyle{D_{k+1}\over D}$. Let $\tilde{D}$ be the gcd of the
$(e,e)$ minors of the matrix
  $[w_1^T,\cdots,w_{e}^T,(i.v)^T]$. If
  $i.v\in G_k$, then $\tilde{D}=D_{k+1}$. But $\tilde{D}=$
  gcd$(id_1,\cdots,id_e,D_{k+1})$, in particular $D_{k+1}$ divides
  gcd$(id_1,\cdots,id_e,iD_{k+1})=i.D$ which is a contradiction
  because $i.D <D_{k+1}$.$\blacksquare$

\end{demostracion}

\vskip0.2cm

\noindent Since $D_1 >\cdots >D_{s+1}$, it follows that $G_0\subset G_1\subset \cdots \subset G_s$. We also have the following:

\begin{proposicion}{\rm i) For all $1\leq k\leq s, e_k$ is the index of
    $G_{k-1}$ in  $G_{k}$.

\vskip0.2cm

ii)  For all $k=1,\cdots,s, e_kv_{e+k}\in G_{k-1}$ and $(e_k-i)v_{e+k}\notin G_{k-1}$ for all $1\leq i <e_k$.

\vskip0.2cm

iii) Given $0\leq k\leq s$ and $v\in {G}_k$, there exist unique
integers
$\lambda_1,\ldots,\lambda_e,\lambda_{e+1},\ldots,\lambda_{e+k}$
such that
$v=\sum_{i=1}^{e+k}\lambda_iv_i$
and $0\leq \lambda_i <e_i$ for all $i=e+1,\ldots,e+k$ (we call this
representation the standard representation with respect to $v_1,\cdots,v_{e+k}$).

\vskip0.2cm

}
\end{proposicion}

\begin{demostracion}{.} i) is obvious and ii) results from Proposition 2.1.
ii). To prove iii), we first prove the existence: let $v=\sum_{i=1}^{e+k}c_iv_i$ where $c_i\in
\mathbb{Z}$ for all $1\leq i\leq e+k$. If $k=0$, then the assertion is clear. Assume
that $k\geq 1$, and that $c_{e+k} <0$. Let $p\in \mathbb{N}^*$ be
such that $0\leq pe_{k}+c_{e+k}< e_{k}$. We have:

$$
v=\sum_{i=1}^{e+k-1}c_iv_i+(c_{e+k}+pe_{k}-pe_{k})v_{e+k}
$$

\noindent since $e_{k}v_{e+k}\in G_{k-1}$, then so is for $-pe_kv_{e+k}$. In particular we can rewrite $v$ as:

$$
v=\sum_{i=1}^{e+k-1}\tilde{c}_iv_i+(\tilde{c}_{e+k})v_{e+k}
$$

\noindent and $0\leq \tilde{c}_{e+k}=pe_{k}+c_{e+k}<e_{k}$. Since
$\sum_{i=1}^{e+k-1}\tilde{c}_iv_i\in G_{k-1}$, then we get the result
by induction on $k$.

\vskip0.1cm

\noindent To prove the uniqueness, let $v=\sum_{i=1}^{e+k}a_iv_i=\sum_{i=1}^{e+k}b_iv_i$
  where for all $i=e+1,\ldots,e+k, 0\leq a_i,b_i <e_i$ and let $j$
  be the greatest integer such that $a_j-b_j\not=0$. Suppose that $j\geq e+1$
  and also that $a_j-b_j> 0$, then

  $$
  (a_j-b_j)v_j=\sum_{i=1}^e(b_i-a_i)+(b_{e+1}-a_{e+1})v_1+\ldots+
  (b_{j-1}-a_{j-1})v_{j-1}\in
  {G}_{j-1}
  $$

  \noindent and $0<a_j-b_j < e_j$. This contradicts ii).
\end{demostracion}$\blacksquare$

\noindent Note that the results of Propositions 2.1. and 2.2. hold assuming
only that the condition (*) of page 2 is satisfied. This will not be the case in the
following Corollary.

\begin{corolario} {\rm Let  $0\leq k\leq s$ and let $v\in {G}_k$. Let

$$
v=\sum_{i=1}^{e+k}\lambda_iv_i
$$

\noindent be the standard representation with respect to
$v_1,\cdots,v_{e+k}$. The vector $v\in \Gamma(v_1,\cdots,v_{e+k})$
if and only if $\lambda_i\geq 0$ for all $i=1,\cdots,e$.}

\end{corolario}

\begin{demostracion}{.} If  $\lambda_i\geq 0$ for all $i=1,\cdots,e$,
  then clearly $v\in \Gamma(v_1,\cdots,v_{e+k})$. Conversely,
  suppose that  $v\in \Gamma(v_1,\cdots,v_{e+k})$, then
  $v=\sum_{i=1}^{e+k}\mu_iv_i$ where $\mu_i\geq 0$ for all
  $i=1,\cdots,e+k$. We shall construct the standard representation of $v$ as in the Proposition above: if $0\leq \mu_i <e_i$ for all $i=e+1,\ldots,e+k$, then it is over. Assume that $\mu_i\geq e_i$ for some $i\geq e+1$ and let $e+j$ be the greatest element with this property. Write $\mu_{j}=pe_j+\tilde{\mu_j}$, where $0\leq \tilde{\mu_j} < e_j$. But $e_jv_{j}\in \Gamma(v_1,\ldots,v_e,v_{e+1},\ldots,v_{j-1})$, in particular $e_jv_{j}=\sum_{i=1}^{j-1}\tilde{\lambda_i}v_{i}$. We finally rewrite $v$ in the following form:

  $$
  v=\sum_{i=1}^{e+k}\tilde{\lambda_i}v_i
  $$

  \noindent where $\tilde{\lambda_i}\geq 0$ and $0\leq \tilde{\lambda_i}<e_i$ for all $i=j,\ldots e+k$. Finally, we get the result by  an easy induction.$\blacksquare$

\end{demostracion}

\section{Proof of Theorem 1. and applications}

\noindent {\bf Proof of Theorem 1.} Let the notations be as in Section 1. and let
$g(\underline{v})=\sum_{k=1}^s(e_k-1)v_{e+k}-\sum_{i=1}^ev_i$. Clearly $g(\underline{v})\in
G(\underline{v})$, and by corollary 2.3., $g(\underline{v})\notin \Gamma(\underline{v})$. Let $u\in C_e-\lbrace (0,\cdots,0)\rbrace$ and let $v=g(\underline{v})+u$. Assume that $v\in G(\underline{v})$ and let

$$
v=\sum_{k=1}^{e+s}\theta_kv_{k}
$$

\noindent be the standard representation of $v$ and recall that $0\leq \theta_{e+k}< e_k$ for all $k=1,\cdots,s$. We have:

$$
\sum_{k=1}^s(e_k-1-\theta_{e+k})v_{e+k}+u=(\theta_1+1)v_1+\cdots+(\theta_e+1)v_e
$$

\noindent But $\sum_{k=1}^s(e_k-1-\theta_{e+k})v_{e+k}+u\in C_e$, which
implies that $\theta_k+1>0$ for all $k=1,\cdots,e$, in particular
$\theta_k\geq 0$ for all $k=1,\cdots,e$, consequently $g(\underline{v})+u\in \Gamma(\underline{v})$.$\blacksquare$


\begin{definicion} {\rm Suppose that $D_{s+1}=1$,
    i.e. $G(\underline{v})=\mathbb{Z}^e$, and let N$(C_e)$ be the
    set of the compact faces of the convex hull of $\bigcup_{w\in C_e}w+C_e$. Let $w_1,\cdots,w_r\in \mathbb{N}^e$ be the
    set of integral vectors of N$(C_e)$. For all $v\in C_e$, there is
    $1\leq k\leq r$ such that $v\in w_k+C_e$. In particular, for all
    $v\in g(\underline{v})+C_e$, if $v\in C_e$, then there exists
    $1\leq k\leq r$ such that $v\in (g+w_k)+C_e$. The set $\lbrace g(\underline{v})+w_1,\cdots,g(\underline{v})+w_r\rbrace$ is called the conductor of $\Gamma(\underline{v})$.}
\end{definicion}

\begin{corolario}{\rm Let $\underline{v}=(v_1,\cdots,v_{e+s})$ be as
    in Section 1. and let
    $A=[v_1^T,\cdots,v_{e+s}^T]$. Consider the Diophantine equation
    $A.X=B$ where $B\in{\mathbb{N}}^e$. By Theorem 1., if
    $B\in g(\underline{v})+C_e$, then $B\in\Gamma(\underline{v})$, in
    particular the Diophantine equation $A.X=B$ has a solution in
    ${\mathbb{N}}^{e+s}$}
\end{corolario}

\subsection{The semigroup of a curve singularity}

\noindent Let ${\bf K}$ be an algebraically closed field of
    characteristic zero and let $f=y^n+a_1(x)y^{n-1}+\cdots+a_n(x)$ be
    a nonzero element of ${\bf K}[[x]][y]$. Suppose that $f$ is
    irreducible. By Newton theorem, there exists
    $y(t)=\sum_{p}c_pt^p\in {\bf K}[[t]]$ such that
    $f(t^n,y(t))=0$. Furthermore, $f(t^n,y)=\prod_{w\in
      U_n}(y-y(wt))$, where $U_n$ denotes the group of roots of unity
    in ${\bf K}$. Given a nonzero polynomial $g\in {\bf K}[[x]][y]$,
    we set int$(f,g)=O_tg(t^n,y(t))$, where $O_t$ denotes the
    $t$-order. The set of int$(f,g),0\not= g\in {\bf K}[[x]][y]$ is a
    numerical semigroup, denoted $\Gamma(f)$. Let $m_0=n=d_1$ and for
    all $k\geq 1$, let $m_k={\rm inf}\lbrace p|c_p\not=0$ and
    $m_{k-1}$ does not divide $p\rbrace$. There exists $h\geq 1$ such
    that $d_{h+1}=1$. The set $\lbrace m_1,\ldots,m_h\rbrace$ is called the
    set of Newton-Puiseux exponents of $f$. With these notations,
    $\Gamma(f)$ is generated by $r_0,r_1,\cdots,r_h$, where $r_0=n, r_1=m_1$
    and for all $2\leq k\leq h$:

$$
r_k= r_{k-1}{d_{k-1}\over d_k}+m_k-m_{k-1}
$$

\noindent and it is well known that
$r_ke_k\in\Gamma(r_0,\ldots,r_{k-1})$ for all $1\leq k\leq h$. Conversely,
let $r_0<r_1,\cdots <r_h$ be a given sequence of relatively prime
nonnegative integers. Let $d_1=r_0$ and for all $1\leq k\leq h$ let $d_{k+1}={\rm
  gcd}(r_{k},d_k)$. If
$r_{k}>r_{k-1}\displaystyle{{d_{k-1}\over d_{k}}}$ for all $1\leq
k\leq h$, then
$r_0,\cdots,r_h$ generate the semigroup of an irreducible element of
${\bf K}[[x]][y]$ (see [11]).

\vskip0.1cm

\noindent Let $f$ be as above, and let $r_0,\cdots,r_h$ be the set of
generators of $\Gamma(f)$. Let $e_k=\displaystyle{d_k\over d_{k+1}}, k=1,\cdots,h$,
and let:

$$
c=\sum_{k=1}^h(e_k-1)r_k-r_0+1
$$

\noindent then $c$ is the conductor of $\Gamma(f)$,
i.e. $g=c-1\notin\Gamma(f)$ and $c+{\mathbb{N}}\subseteq
\Gamma(f)$. The ideal $(t^c)$ is the conductor ideal of ${\bf
  K}[[x]][y]/(f)$ into its intergal closure ${\bf K}[[t]]$, and $c$ is
also the Milnor number of $f$, i.e. $c={\rm rank}_{\bf K}{\bf
  K}[[x,y]]/(f_x,f_y)$, where $f_x$ (resp. $f_y$) denotes the
$x$-derivative (resp. the $y$-derivative) of $f$. Furthermore,
the cardinality of $\mathbb{N}-\Gamma(f)$ (the set of gaps of
$\Gamma(f)$)  is $\displaystyle{c\over 2}$.

\subsection{The semigroup of a quasi-ordinary polynomial}

\noindent Let ${\bf K}$ be an algebraically closed field of
    characteristic zero and let
    $f=y^n+a_1(\underline{x})y^{n-1}+\cdots+a_n(\underline{x})$ be an irreducible
    element of ${\bf K}[[\underline{x}]][y]={\bf
      K}[[x_1,\cdots,x_e]][y]$ and assume that the discriminant $D_y(f)$
    of $f$, defined to be the $y$-resultant of $f$ and its
    $y$-derivative $f_y$, is of the form
    $x_1^{N_1}\cdots x_e^{N_e}(a+u(\underline{x}))$, where
    $N_1,\cdots,N_e\in\mathbb{N}, a\in {\bf K}^*$, and
    $u(\underline{0})=0$ (such a polynomial is called a quasi-ordinary
    polynomial). By [1], there exists
    $y(\underline{t})=y(t_1,\cdots,t_e)=\sum_{p\in {\mathbb N}^e}c_p{\underline{t}}^p\in {\bf K}[[t_1,\cdots,t_e]]$
    such that $f(t_1^n,\cdots,t_e^n,y(\underline{t}))=0$. Furthermore,
there exist $n$ distinct elements $(\underline{w}^1,\cdots,\underline{w}^n)\in U_n^e$, where
$U_n$ denotes the group of roots of unity in ${\bf K}$, such that
$f(t_1^n,\cdots,t_e^n,y)=\prod_{k=1}^n(y-y(w_1^kt_1,\cdots,w_e^kt_e))$.

\noindent Given a nonzero element $g\in{\bf K}[[\underline{x}]][y]$,
  we define $O(f,g)$ to be the maximal element with respect to the
  lexicographical order of
  the initial form of $g(t_1^n,\cdots,t_e^n,y(\underline{t}))$. The set of $O(f,g),0\not=
  g\in {\bf K}[[\underline{x}]][y]$, is a semigroup of $\mathbb{N}^e$,
  denotes $\Gamma(f)$. Let
Supp$(y(\underline{t}))=\lbrace p| c_p\not=0\rbrace$. In
[7], J. Lipman proved the existence of $m_1,\cdots,m_h\in$
Supp$(y(\underline{t}))$ such that the following hold:

i) $m_1 < m_2< \cdots <m_h$, where $<$ means $<$ coordinate-wise.

ii) Let $M_0=(n\mathbb{Z})^e$ and for all $k=1,\cdots,h$, let
$M_k=M_0+\sum_{i=1}^km_i\mathbb{Z}$. We have $M_0\subset
M_1\subset\cdots \subset M_h$. Furthermore, for all $p\in{\rm
  Supp}(y(\underline{t})), p\in \sum_{p\in
  m_k+\mathbb{N}^e}m_k\mathbb{Z}$.

\noindent Let $D_1=n^e$ and for all $k=1,\cdots,h$, let $D_{k+1}$ be
the gcd of $e\times e$ minors of the matrix $[nI_e,m_1^T,\cdots,m_k^T]$, where $I_e$ denotes the $(e,e)$ unit matrix. By conditions i), ii), we have $D_1 >\cdots
>D_{h+1}$, furthermore $D_{h+1}=n^{e-1}$ (see [3]). Let
$r_0^1,\cdots,r_0^e$ be the canonical basis of $(n\mathbb{Z})^e$, and
define the sequence $r_1,\ldots,r_h$ by $r_1=m_1$ and for all
$k=2,\cdots,h$:

$$
r_k=r_{k-1}{D_{k-1}\over D_k}+(m_k-m_{k-1})
$$

\noindent then $r_0^1,\cdots,r_0^e,r_1,\cdots,r_h$ generate
$\Gamma(f)$ and $\displaystyle{r_k{D_k\over D_{k+1}}}\in
\Gamma(r_0^1,\ldots,r_0^e,r_1,\ldots,r_{k-1})$ for all $1\leq k\leq h$. Furthermore, for all $k=1,\cdots,h$, if $\tilde{D}_{k+1}$
denotes the gcd of the $(e,e)$ minors of the matrix $[nI_e,
r_1^T,\cdots,r_{k}^T]$, then $\tilde{D}_{k+1}=D_{k+1}$.

\medskip

\noindent Note that in this situation, the convex cone generated by
$r_0^1,\cdots,r_0^e$ is nothing but $\mathbb{R}_+^e$, and
$C_e=(\mathbb{R}_+^*)^e$.

\medskip

\noindent Set $e_k=\displaystyle{{D_k\over D_{k+1}}}$ for all $k=1,\cdots,h$. The Frobenius vector of $\Gamma(f)$ is:

$$
g=\sum_{k=1}^h (e_k-1)r_k-\sum_{k=1}^e r_0^k
$$

\noindent in particular, for all $u\in (\mathbb{R}_+^*)^e$, if $g+u
\in G(r_0^1,\cdots,r_0^e,r_1,\cdots,r_h)$, then $g+u\in \Gamma(f)$.

\subsection{Numerical examples}

\begin{exemple}{\rm (See Fig. 2) Let $\underline{v}=(v_1,v_2,v_3)=(4,6,7)$. With the notations of Section 1. we have $D_1=4,D_2=2,D_3=1,
\displaystyle{e_1={D_1\over D_2}=2, e_2={D_2\over D_3}=2}$. The Frobenius vector  of $\Gamma(\underline{v})$ is:

$$
g(\underline{v})=(e_1-1)v_2+(e_2-1)v_3-v_1= 9
$$

\noindent and the conductor of $\Gamma(\underline{v})$ is $c=g(\underline{v})+1=10$. Note that $\mathbb{N}-\Gamma(\underline{v})=\lbrace 1,2,3,5,9\rbrace$, whose cardinality is $\displaystyle{{c\over 2}=5}$.}
\end{exemple}

\begin{exemple}{\rm (See Fig. 3) Let
    $\underline{v}=(v_1,v_2,v_3,v_4)=((8,0),(0,8),(2,2),(12,8))$. With
    the notations of Section 1. we have $D_1=64$, $D_2$-the gcd of
    the $(2,2)$ minors of the matrix $[8I_2,(2,2)^T]$- is $16$, and $D_3$-the gcd of the $(2,2)$ minors of the matrix $[8I_2,(2,2)^T,(12,8)^T]$- is  $8$. Finally,
$\displaystyle{e_1={D_1\over D_2}=4, e_2={D_2\over D_3}=2}$. The Frobenius vector  of $\Gamma(\underline{v})$ is:

$$
g(\underline{v})=(e_1-1)v_3+(e_2-1)v_4-v_1-v_2= 3v_3+v_4-v_1-v_2=(10,6)
$$

\noindent Let $u=v_2=(2,2)$, then
$g+u=g+v_2=(e_1-1)v_3+(e_2-1)v_4-v_1=(10,14)\notin
\Gamma(\underline{v})$. In fact, $u$ belongs to a cell of
cone$(v_1,v_2)$ of dimension $1$.}

\end{exemple}

\begin{exemple}{\rm (See Fig. 4) Let
    $\underline{v}=(v_1,v_2,v_3,v_4)=((4,6),(6,3),(8,10),(3,4))$. With
    the notations of Section 1. we have $D_1=24$, $D_2$-the gcd of the
    $(2,2)$ minors of the matrix $[(4,6)^T,(6,3)^T,(8,10)^T]$- is $4$,
    and $D_3$-the gcd of the $(2,2)$ minors of the matrix $[(4,6)^T,(6,3)^T,(8,10)^T,(3,4)^T]$- is  $1$. Finally,
$\displaystyle{e_1={D_1\over D_2}=6, e_2={D_2\over D_3}=4}$. The Frobenius vector  of $\Gamma(\underline{v})$ is:

$$
g(\underline{v})=(e_1-1)v_3+(e_2-1)v_4-v_1-v_2= (49,52)-(10,9)=(39,53)
$$

\noindent In this example, since $D_3=1$, then
$G(\underline{v})={\mathbb{Z}}^2$. In particular, for all $v\in g(\underline{v})+C_e, v\in\Gamma(\underline{v})$. Furthermore, for all $v\in C_e$,
 $v\in (1,1)+C_e$, hence the conductor of $\Gamma(\underline{v})$ is $g(\underline{v})+(1,1)=(40,54)$.
}

\end{exemple}

\begin{exemple}{\rm (See Fig. 5) Let $\underline{v}=(v_1,v_2,v_3)=((1,3),(3,2),(1,1))$. With the notations of Section 1. we have $D_1=7$, $D_2$-the gcd of the $(2,2)$ minors of the matrix $[(1,3)^T,(3,2)^T,(1,1)^T]$- is $1$. Finally,
$\displaystyle{e_1={D_1\over D_2}}=7$. The Frobenius number  of $\Gamma(\underline{v})$ is:

$$
g(\underline{v})=(e_1-1)v_3-v_1-v_2= (6,6)-(4,5)=(2,1)
$$

\noindent In this example, since $D_2=1$, then
$G(\underline{v})={\mathbb{Z}}^2$. In particular, for all $v\in g(\underline{v})+C_e, v\in\Gamma(\underline{v})$. Furthermore, for all $v\in C_e$,
 $v\in [(1,1)+C_e]\cup [(1,2)+C_e]$, hence the conductor of $\Gamma(\underline{v})$ is $\lbrace g(\underline{v})+(1,1), g(\underline{v})+(1,2)]\rbrace=\lbrace (3,2),(3,3)\rbrace$.
}

\end{exemple}

\bigskip
\bigskip
\bigskip
\bigskip

\unitlength=0.5cm
\begin{picture}(6,4)
\put(0,6){\vector(1,0){16}} 
\put(0,6){\circle*{.2}}\put(0,6.5){\shortstack{0}}
\put(1,6){\circle{.2}}
\put(2,6){\circle{.2}}
\put(3,6){\circle{.2}}
\put(4,6){\circle*{.2}}\put(3.9,6.5){\shortstack{4}}
\put(5,6){\circle{.2}}
\put(6,6){\circle*{.2}}\put(5.9,6.5){\shortstack{6}}
\put(7,6){\circle*{.2}}\put(6.9,6.5){\shortstack{7}}
\put(8,6){\circle*{.2}}\put(7.9,6.5){\shortstack{8}}
\put(9,6){\circle{.2}}\put(8.9,6.5){\shortstack{$g$}}
\put(10,6){\circle*{.2}}\put(9.8,6.5){\shortstack{10}}
\put(11,6){\circle*{.2}}
\put(12,6){\circle*{.2}}
\put(13,6){\circle*{.2}}
\put(14,6){\circle*{.2}}
\put(18,6){\shortstack{Fig. 2}}
\end{picture}

\unitlength=0.5cm
\begin{picture}(6,4)

\put(0,-1.5){\vector(0,1){10}} 
\put(-0.5,-1){\vector(1,0){10}}

\put(6,-0.7){\shortstack{$v_1$}}
\put(0.2,4.8){\shortstack{$v_2$}}
\put(7,3.2){\circle*{.3}}\put(7,2.3){\shortstack{$g$}}
\put(7,3.2){\vector(1,0){4}}
\put(7,3.2){\vector(0,1){4}}
\put(0,-1){\vector(1,0){6}}
\put(0,-1){\vector(0,1){5.8}}
\put(0,-1){\vector(1,1){2}}\put(2,1.2){\shortstack{$v_3$}}
\put(0,-1){\vector(3,2){8}}\put(8,4.2){\shortstack{$v_4$}}
\put(12,2.3){\shortstack{Fig. 3}}
\end{picture}

\unitlength=0.5cm
\begin{picture}(6,4)
\put(17,2.5){\vector(0,1){10}} 
\put(16.5,3){\vector(1,0){11}}
\put(17,3){\vector(2,3){2}}\put(18.2,6.3){\shortstack{$v_1$}}
\put(17,3){\vector(22,10){2}}\put(20.2,8.7){\shortstack{$v_3$}}
\put(17,3){\vector(4,5){5}}\put(18.5,4.5){\shortstack{$v_4$}}
\put(17,3){\vector(3,4){1.5}}\put(18.3,3.5){\shortstack{$v_2$}}
\put(23.9,9.2){\shortstack{$\small{g}$}}

\put(24,10){\circle*{.3}}
\put(24.5,10.5){\circle{.2}}
\put(24,10){\vector(2,3){3}}
\put(24.5,10.5){\line(2,3){2.7}}
\put(24,10){\vector(2,1){3}}
\put(24.5,10.5){\line(2,1){3}}
\put(28,5.5){\shortstack{Fig. 4}}
\end{picture}

\unitlength=0.5cm
\begin{picture}(6,4)
\put(0,-3){\circle*{.3}}
\put(0,-3.5){\vector(0,1){10}} 
\put(-0.5,-3){\vector(1,0){10}}
\put(0,-3){\vector(1,3){2}}
\put(2,3){\shortstack{$v_1$}}
\put(0,-3){\vector(3,2){5}}
\put(5,0.8){\shortstack{$v_2$}}
\put(0,-3){\vector(1,1){5}}
\put(5,2){\shortstack{$v_3$}}
\put(1,-2){\circle*{0.2}}
\put(1,-1){\circle*{0.2}}
\put(10,0.8){\shortstack{Fig. 5}}
\put(2,-2.5){\shortstack{$\small{g}$}}\put(2,-2){\circle*{0.2}}
\put(2,-2){\vector(1,3){2}}
\put(2,-2){\vector(3,2){3}}
\put(3,-1){\circle{0.2}}
\put(3,-1){\line(1,3){2}}
\put(3,-1){\line(3,2){2.5}}
\put(3,0){\circle{0.2}}
\put(3,0){\line(1,3){2}}
\put(3,0){\line(3,2){3}}
\put(3,-1){\line(0,1){1}}

\end{picture}

\vskip2cm

\begin{nota}{\rm i) Let $\underline{v}=(v_1,\cdots,v_{e+s})$ et let the notations be as in Section 1. and suppose
that $D_1\geq D_2\geq\cdots \geq D_{s+1}$. Let $(v_{e+i_1},\cdots,v_{e+i_t})$ be the maximal set of $\lbrace v_{e+1},\cdots,v_{e+s}\rbrace$ such that for all $1\leq k\leq t, v_{e+i_k}\notin G(v_1,\cdots,v_e,\cdots,v_{e+i_k-1})$ and let
$g=\sum_{k=1}^t(e_{i_k}-1)v_{e+i_k}-\sum_{k=1}^ev_k$. For all $v\in g+C_e$, if $v\in G(\underline{v})=G(v_1,\cdots,v_e,v_{i_1},\cdots,v_{i_t})$, then $v\in \Gamma(v_1,\cdots,v_e,v_{i_1},\cdots,v_{i_t})\subseteq \Gamma(\underline{v})$. In general, the vector $g$ need not to be the "smallest" one with this property. Let for example $\underline{v}=(4,6,7,9)$: for all $v\geq 9, v\in \Gamma(\underline{v})$, but the Frobenius number of $\Gamma(\underline{v})$ is $5$. General subsemigroups of $\mathbb{N}$ and their Frobenius numbers have been studied by many authors (see [8] and references).

ii) If $\underline{v}=(8,10,11)$, then $g(\underline{v})=(e_1-1)v_2+(e_-1)v_3-v_1=30+11-8=33=3.11\in \Gamma(8,10,11)$. In this example, condition (**) of page 2 is not satisfied, since $e_2v_3=22\notin\Gamma(8,10)$. }
\end{nota}

\end{document}